\documentclass[10pt]{article}
\usepackage{amsmath,amssymb,amsthm}
\usepackage{hyperref}
\usepackage{units}
\usepackage{color}
\usepackage[T1]{fontenc}
\usepackage[utf8]{inputenc}
\usepackage{authblk}
\usepackage{bm}

\begin{document}
\title{The Standard Representation of the Symmetric Group $S_n$ over the Ring of Integers\thanks{%
Mathematics Subject Classifications: 20C30\newline Keywords: standard representation, symmetric group, casimir invariant }}

\author[]{Kunle Adegoke\thanks{Corresponding author: adegoke00@gmail.com}}
\author[2]{Olawanle Layeni}
\author[]{Rauf Giwa}
\author[]{Gbenga Olunloyo}
\affil{Department of Physics and Engineering Physics, Obafemi Awolowo University, Ile-Ife, Nigeria}
\affil[2]{Department of Mathematics, Obafemi Awolowo University, Ile-Ife, Nigeria}
\renewcommand\Authands{ and }
\date{}

\maketitle
\tableofcontents
\begin{abstract}
\noindent In this paper we give a Casimir Invariant for the Symmetric group $S_n$. Furthermore we obtain and present, for the first time in the literature, explicit formulas for the matrices of the standard representation in terms of the matrices of the permutation representation.
\end{abstract}

\newtheorem{thm}{Theorem}
\newtheorem{cor}[thm]{Corollary}
\newtheorem{lem}[thm]{Lemma}
\newtheorem{prop}[thm]{Proposition}
\newtheorem{ax}{Axiom}

\newtheorem{defn}{Definition}[section]

\newtheorem{rem}{Remark}[section]

\section{Introduction}

The symmetric group $S_n$, whose elements consist of the set of all permutations on $n$ symbols is of central importance to mathematics and physics \cite{hamermesh}. Cayley's theorem states that every group is isomorphic to a subgroup of the symmetric group on that group. In physics the classification of atomic and nuclear states depends essentially on the properties of $S_n$ \cite{hamermesh}. The representation theory of the symmetric group is a well studied subject~\cite{hamermesh, burnside, fulton91}. The partitions of $n$ or equivalently Young diagrams of size $n$ are the natural ways in which to parametrize the irreducible representations of $S_n$ \cite{wikisymmetric1}. This paper is concerned, not with the general irreducible representations of $S_n$, but, more specifically, with the so-called standard representation of the symmetric group, formally obtained from the \mbox{$n - 1$} dimensional subspace of vectors whose sum of coordinates is zero in the basis set of a permutation representation. The path taken in this work shall however be non-group theoretic. For example, we will not be concerned with Young diagrams. 

The standard representation of $S_n$ is important for the following reason: For $n \ge 7$, the permutation representation, the trivial (identical) representation, the sign representation, the standard representation and another $n-1$ dimensional irreducible representation found by tensoring with the sign representation are the only lowest-dimensional irreducible representations of $S_n$ \cite{wikisymmetric1}. All other irreducible representations have dimension at least $n$. While it is a fact that all irreducible representations of $S_n$ can be found, using Frobenius formula (\cite{hamermesh}, pp 189), for example, there are no known explicit formulas for the standard representation. The main result of this paper is the derivation of such formulas, which now make it possible to write down the standard representation matrices directly from those of the permutation representation.

\section{The Permutation Representation}

Denote the $n!$ elements of $S_n$ by $A^k$, $k = 1,2, \ldots, n!$, such that, in usual notation,

\begin{equation}
A^k  = \left( {\begin{array}{*{20}c}
   1 & 2 & 3 &  \cdots  &  \cdots  &  \cdots  & n  \\
   {a_1^k } & {a_2^k } & {a_3^k } &  \cdots  &  \cdots  &  \cdots  & {a_n^k }  \\
\end{array}} \right)\,,
\end{equation}
and $1\le a_i^k\le n$, all $a_i^k$ being distinct. For simplicity, and since no ambiguity can result, we will use the same symbol $A^k$ for the representation matrices. Then, in the permutation representation, the $n\times n$ matrices $A^k$ are given, through their elements, by:

\begin{equation}
(A^k )_{ij}  = \delta _{a_i^k ,j} ,\quad i = 1 \ldots n,\;j = 1 \ldots n\,,
\end{equation}

which is clearly a unitary representation, since 

\[
\left( {(A^k )^ \dagger  A^k } \right)_{ij}  = \sum\limits_{r = 1}^n {\left( {(A^k )^ \dagger  } \right)_{ir} A_{rj}^k }  = \sum\limits_{r = 1}^n {A_{ri}^k A_{rj}^k }  = \sum\limits_{r = 1}^n {\delta _{a_r^k ,i} \delta _{a_r^k ,j} }  = \delta _{ij}\,.
\]

\subsection{A Casimir Invariant for $S_n$ in the permutation representation}

\begin{thm}
The $n\times n$ matrix $C$, with elements $C_{ij}=1-\delta_{ij}$ is a Casimir Invariant of $S_n$ in the permutation representation.
\end{thm}

\textbf{Proof} We require to prove that $C$ commutes with every $A^k$.\\
\begin{equation}
\begin{split}
(CA^k )_{ij}  = \sum\limits_{r = 1}^n {C_{ir} (A^k )_{rj} }  &= \sum\limits_{r = 1}^n {(1 - \delta _{ir} )(A^k )_{rj} }\\
&= \sum\limits_{r = 1}^n {(A^k )_{rj} }  - \sum\limits_{r = 1}^n {\delta _{ir} (A^k )_{rj} }\\
&= \sum\limits_{r = 1}^n {(A^k )_{rj} }  - (A^k )_{ij}\\
&= 1 - (A^k )_{ij}\,.
\end{split}
\end{equation}

A similar calculation gives,

\begin{equation}
(A^k C)_{ij} = \sum\limits_{r = 1}^n {(A^k )_{ir} C_{rj} }  = \sum\limits_{r = 1}^n {(A^k )_{ir} }  - (A^k )_{ij}= 1 - (A^k )_{ij}\,.
\end{equation}

We see, therefore, that $A^kC=CA^k$, so that $C$ is a Casimir Invariant of the symmetric group.


\section{The Standard Representation}

Since the Casimir invariant $C$, obtained in the previous section, is not proportional to the identity, Schur's lemma tells us that the permutation representation is not irreducible, a well-known fact. It therefore remains to find the matrix $P$ which diagonalizes $C$. First we prove a lemma.

\textbf{Lemma 1} {\em The nonsingular $n\times n$ matrix $P$ with elements \mbox{$P_{ij}=\delta_{j1}-(1-\delta_{j1})\delta_{i1}+\delta_{i,n+2-j}$} has the inverse $P^{-1}$ where \mbox{$n(P^{-1})_{ij}=2\delta_{i1}-1+n\delta_{i,n+2-j}$}}

\textbf{Proof}
\begin{equation}
\begin{split}
n(PP^{ - 1} )_{ij}&= \sum\limits_{r = 1}^n {P_{ir} n(P^{ - 1} )_{rj} }\\  
&= \sum\limits_{r = 1}^n {(\delta _{r1}  - \delta _{i1}  + \delta _{r1} \delta _{i1}  + \delta _{i,n + 2 - r} )}\\
&\qquad\qquad\qquad\times (2\delta _{r1}  - 1 + n\delta _{r,n + 2 - j} )
\end{split}
\end{equation}

It is straightforward to write out the terms and evaluate the summation termwise. One merely needs to note that 

\[\sum\nolimits_{r = 1}^n {\delta _{r,n + 2 - j} }  = 1 - \delta _{j1}\] and \[\sum\nolimits_{r = 1}^n {\delta _{i,n + 2 - r} } \delta _{r,n + 2 - j}  = \delta _{ij} - \delta _{i1} \delta _{j1}\,.\] One then finds \mbox{$(PP^{-1})_{ij}=\delta_{ij}$}, which establishes the claim.

\subsection{Diagonal form of the Casimir Invariant}

\begin{thm}
The matrix $P$ given in Lemma 1, diagonalizes the Casimir Invariant, $C$.
\end{thm}
\textbf{Proof.} We wish to compute
\[
D=P^{-1}CP\,.
\]
Now
\[
D_{ij}  = \sum\limits_{r = 1}^n {\sum\limits_{s = 1}^n {(P^{ - 1} )_{ir} C_{rs} P_{sj} } }\,. 
\]

Substituting the matrix elements, expanding and evaluating the sums, we find after some algebra, that

\begin{equation}
D_{ij}=n\delta _{i1} \delta _{j1}-\delta _{ij}\,.
\end{equation}

Thus we see that $D$ is a diagonal matrix, as claimed, with the entry `$n-1$' in row $1$, column $1$ and the remaining diagonal elements being $-1$.

The matrix $P$, above, which diagonalizes $C$ will block-diagonalize the matrices $A^k$.

\subsection{Similarity Transformation of $A^k$: the standard representation}

Using the matrix elements of $P$, $P^{-1}$ and $A^k$, it is not difficult to obtain the interesting result:

\begin{equation}\label{equ.m5z7nvb}
\begin{split}
(P^{ - 1} A^k P)_{ij} & = \sum\limits_{r = 1}^n {\sum\limits_{s = 1}^n {(P^{ - 1} )_{ir} A_{rs}^k P_{sj} } }\\
&= \delta _{j1} \delta _{i1}  + (1 - \delta _{i1} )(1 - \delta _{j1} )[A_{n + 2 - i,n + 2 - j}^k  - A^k_{n + 2 - i,1}]
\end{split}
\end{equation}

We see from \eqref{equ.m5z7nvb} that each matrix $P^{ - 1} A^k P$ is block diagonal, being the direct sum of a $1\times 1$ matrix with entry $1$ and an $(n-1)\times (n-1)$ matrix $B^k$, with elements

\begin{equation}
B_{ij}^k  = A_{n + 1 - i,n + 1 - j}^k  - A_{n + 1 - i,1}^k \;,\quad i= 1, \ldots ,n - 1,\;j = 1, \ldots ,n - 1\,.
\end{equation}

The $1\times 1$ matrices correspond to the identical (trivial) representation in which every element of $S_n$ is sent to the one-dimensional identity matrix, while the $B^k$ matrices correspond to the irreducible $n-1$ dimensional standard representation.

\section{Conclusion}

In this paper we have shown that the operator $C$ with matrix elements \mbox{$C_{ij}=1-\delta_{ij}$} is a Casimir Invariant for the symmetric group $S_n$. We also showed that if $A^k,\;k=1,2,\ldots, n!$ are the representation matrices for the elements of $S_n$ in the permutation representation, then the matrices $B^k$ for the standard representation of $S_n$ are given by

\[
B_{ij}^k  = A_{n + 1 - i,n + 1 - j}^k  - A_{n + 1 - i,1}^k \;,\quad i= 1, \ldots ,n - 1,\;j = 1, \ldots ,n - 1\,.
\]

\end{document}